\documentclass[11pt]{article} 

\usepackage[margin=1in]{geometry}
\usepackage{setspace}

\normalsize
\usepackage{amsfonts}
\usepackage{amssymb}
\usepackage{amscd}
\usepackage{graphics}
\usepackage[dvips]{graphicx}
\usepackage{epsfig}
\usepackage[cmex10]{amsmath}
\usepackage{array}
\usepackage{eqparbox}
\usepackage{hyperref}
\usepackage{fancyhdr}
\usepackage{appendix}
\usepackage{titling}
\usepackage{enumerate}
\usepackage{tabularx}
\usepackage[super]{nth}
\usepackage{caption} 
\usepackage{subcaption}

\usepackage{booktabs}
\usepackage{float} 

\usepackage[round]{natbib}  
\DeclareMathAlphabet{\mathbsf}{OT1}{cmss}{bx}{n}
\DeclareMathAlphabet{\mathssf}{OT1}{cmss}{m}{sl}
\DeclareMathAlphabet{\mathcsf}{OT1}{cmss}{sbc}{n}

\usepackage{stackengine}
\newcommand{\ie}{{\em i.e.}}
\newcommand{\etc}{{\em etc}}
\newcommand{\eg}{{\em e.g.}}

\newcommand{\apriori}{{\em a priori}}

\newcommand{\secref}[1]{Section~\ref{#1}}
\newcommand{\figref}[1]{Fig.~\ref{#1}}
\newcommand{\tabref}[1]{Table~\ref{#1}}

\newcommand{\keywords}[1]{\textbf{Keywords:} #1}

\makeatletter
\def\blfootnote{\xdef\@thefnmark{}\@footnotetext}
\makeatother

\newcommand{\qed}{\nobreak \ifvmode \relax \else
      \ifdim\lastskip<1.5em \hskip-\lastskip
      \hskip1.5em plus0em minus0.5em \fi \nobreak
      \vrule height0.75em width0.5em depth0.25em\fi}

\usepackage{newlfont}
\date{}
\begin{document}
\title{\textbf{A unified approach to the Behrens--Fisher problem}}
\author{Nagananda K G and Jong Sung Kim\thanks{The authors are with Fariborz Maseeh Department of Mathematics and Statistics, Portland State University, Portland, OR 97201, USA. E-mail: \texttt{\{nanda, jong\}@pdx.edu}.}}
\setlength{\droptitle}{-1.in}
\maketitle
\vspace{-2cm}

\begin{abstract}
A unified framework is presented to study the two-sample Behrens--Fisher problem---testing equality of means when two normal populations have unequal, unknown variances---and a compact expression is derived for the null distribution of the classical test statistic.  Our new approach involves a Mellin--Barnes factorization that decouples the square root of a weighted sum of independent chi-square variates, thereby collapsing a challenging two-dimensional integral to a tractable single-contour integral.  Closing the contour yields a residue series that terminates whenever either sample's degrees of freedom is odd.  A complementary Euler--Beta reduction identifies the density as a Gauss hypergeometric function with explicit parameters, yielding a numerically stable form that recovers Student's $t$ under equal variances.  Ramanujan's master theorem supplies exact inverse-power tail coefficients, which bound Lugannani--Rice saddle-point approximation errors and support reliable tail analyses.  The proposed framework reveals why hypergeometric structure appears, why certain finite-sum cases arise, and how one can pass from the bulk of the distribution to its tails without altering the analytic framework.  Finally, it lets us tabulate exact two-sided critical values over a broad grid of sample sizes and variance ratios that reveal the parameter surface on which the well-known Welch's approximation switches from conservative to liberal, quantifying its maximum size distortion.
\end{abstract}

\section{Introduction}\label{sec:introduction}
Consider independent samples $X_{1},\dots ,X_{n_{1}}\stackrel{\text{iid}}\sim\mathcal N(\mu_{1},\sigma_{1}^{2})$ and $Y_{1},\dots ,Y_{n_{2}}\stackrel{\text{iid}}\sim\mathcal N(\mu_{2},\sigma_{2}^{2})$ from two normal populations.  Let $\nu_{i}=n_{i}-1$, $\bar{X} = \tfrac1{n_1}\sum_{i=1}^{n_1} X_{i}$, $\bar{Y} = \tfrac1{n_2}\sum_{i=1}^{n_2} Y_{i}$, $S_{1}^{2} = \tfrac1{\nu_{1}}\sum_{i=1}^{n_1} (X_{i} - \bar{X})^{2}$, $S_{2}^{2} = \tfrac1{\nu_{2}}\sum_{i=1}^{n_2} (Y_{i}-\bar{Y})^{2}$ denote the degrees of freedom, the sample means and the sample variances, respectively, of the two samples. Also, let $g = \sigma_{1}^{2}/n_{1}$, $h =\sigma_{2}^{2}/n_{2}$ and variance ratio $r = \sigma_1^2/\sigma_2^2$.  Following the pioneering work of \cite{Behrens1929} and \cite{Fisher1935}, the problem of testing $H_{0}\!: \mu_{1} = \mu_{2}$ with $(\sigma_1^2, \sigma_2^2)$ unknown and $\sigma_1^2 \neq \sigma_2^2$ remains an enduring challenge in mathematical statistics.  Several approximations, modifications and practical procedures have been proposed to address this long-standing problem\textemdash for instance, see the work of \cite{Scheffe1970,Bozdogan1986,Best1987,Asiribo1989,Duong1996,Kim1998a,Vangel2005,Dudewicz2007,Chang2008,Nadarajah2017,Chaturvedi2019,Wang2022,Chen2022,Chen2023}.  Other variants  including nonparametric versions of the test (see \cite{Brunner2002,Larocque2010,Konietschke2012,Konietschke2012a,He2022}) and high-dimensional settings (see, for example, \cite{Zhou2017,Cao2019,Zhang2021,Pei2026}) have also been reported in the literature.  Despite these useful and impactful contributions, to the best of our knowledge, no explicit canonical expression for the density function of the Behrens--Fisher statistic 
\begin{eqnarray}  
T = \frac{\bar X-\bar Y}{\sqrt{\frac{S_{1}^{2}}{n_{1}} + \frac{S_{2}^{2}}{n_{2}}}}
\end{eqnarray}
has been established in terms of standard probability distribution functions.  The absence of such a representation complicates exact significance evaluation, power calculations, and asymptotic error control\textemdash particularly in the small-sample, high-variance-imbalance regimes common in modern experimental design.  In what follows, we first explain why the Behrens--Fisher problem remains a persistent difficulty in mathematical statistics, and then in \secref{sec:our_contribution} outline how our work advances its analysis.

Under $H_0$, $\bar{X} - \bar{Y} \sim \mathcal{N}\left(0, g + h\right)$; $\frac{\nu_1S_1^2}{\sigma_1^2} \sim \chi^2_{\nu_1}$, $\frac{\nu_2S_2^2}{\sigma_2^2} \sim \chi^2_{\nu_2}$, and these are independent of $\bar{X}$ and $\bar{Y}$.  The key challenge is that, when forming the ratio for $T$ the denominator involves $\sqrt{S_1^2/n_1 + S_2^2/n_2}$, which is a nonlinear function of two independent $\chi^2$ variables scaled by $\sigma_1^2$ and $\sigma_2^2$ that are unknown and $\sigma_1^2 \neq \sigma_2^2$.  In attempting a closed-form expression for the distribution of \(T\), conditioned on {\apriori} knowledge of \((\sigma_1^2, \sigma_2^2)\), we can write 
\begin{eqnarray*} 
T &=& \frac{\bar{X} - \bar{Y}}{\sqrt{\frac{gW_1}{\nu_1} + \frac{hW_2}{\nu_2}}}, 
\end{eqnarray*} 
where $W_i = \nu_i S_i^2/\sigma_i^2$, $i = 1, 2$.  Under $H_0$, we have $W_i \sim \chi^2_{\nu_i }$, $i = 1, 2$.  This yields 
\begin{eqnarray*} 
T \left| \sigma_1^2, \sigma_2^2 \right. \sim \frac{\mathcal{N}\left(0, g + h\right)}{\sqrt{\frac{g W_1}{\nu_1} + \frac{h W_2}{\nu_2}}}.
\end{eqnarray*} 
When $\sigma_1^2 = \sigma_2^2 = \sigma^2$, and $\sigma^2$ unknown, we have the likelihood ratio statistic 
\begin{eqnarray*} 
T\left| \sigma^2\right. \sim \frac{\left(\bar{X} - \bar{Y}\right)}{\sqrt{S^2_p \left(\frac{1}{n_1} + \frac{1}{n_2}\right)}}, 
\end{eqnarray*} 
where the pooled variance estimate 
\begin{eqnarray*} 
S^2_p = \frac{1}{\nu_1 + \nu_2}\left(\sum_{j=1}^{n_1}(X_j - \bar{X})^2 +  \sum_{j=1}^{n_2}(Y_j - \bar{Y})^2\right).
\end{eqnarray*} 
This implies $T\left| \sigma^2\right. \sim t_{\nu_1 + \nu_2}$, the Student's two-sample $t-$test following a $t-$distribution with $(\nu_1 + \nu_2)$ degrees of freedom (see \cite[Exercise 8.41]{Casella2001}).  However, when $\sigma_1^2 \neq \sigma_2^2$ and $(\sigma_1^2, \sigma_2^2)$ unknown, no such simplification\textemdash including using the pooled variance estimate $S^2_p$\textemdash is possible (see \cite[Exercise 8.42]{Casella2001}).  The distribution of $T = \frac{\bar{X} - \bar{Y}}{\sqrt{\frac{gW_1}{\nu_1} + \frac{hW_2}{\nu_2}}}$ would require integrating over the joint density of \((W_1, W_2)\): 
\begin{eqnarray*} 
f_T(t) = \int_0^\infty \int_0^\infty f_{T|W_1,W_2}(t) f_{W_1}(w_1) f_{W_2}(w_2) \,dw_1\,dw_2,
\end{eqnarray*} 
where each of \(f_{W_1}(w_1)\) and \(f_{W_2}(w_2)\) is the pdf of a $\chi^2$ distribution. The conditional density \(f_{T|W_1,W_2}(t)\) is normal with a variance that depends on \(w_1\) and \(w_2\). Thus, we have 
\begin{eqnarray}
f_T(t) =\!\!\! \int\limits_0^\infty \int\limits_0^\infty \frac{\sqrt{gw_1/\nu_1 + hw_2/\nu_2}}{\sqrt{2\pi(g + h)}} \exp\left(-\frac{t^2\left(gw_1/\nu_1 + hw_2/\nu_2\right)}{2(g + h)}\right) f_{W_1}(w_1) f_{W_2}(w_2)\,dw_1\,dw_2.
\label{eq:behrens_fisher_integral}
\end{eqnarray}
Efforts to solve \eqref{eq:behrens_fisher_integral} does not yield a closed-form expression in terms of elementary distributions. This complexity is the centerpiece of the Behrens--Fisher problem: lifting the equal-variance assumption transforms the otherwise straightforward two-sample $t-$test into a setting where no exact, closed-form distribution is available in elementary terms.  In addressing this problem, \cite{Welch1947} employed the approximation developed by \cite{Satterthwaite1946} to replace the unknown degrees of freedom with a carefully crafted estimate, 
\begin{eqnarray*}
\hat{\nu} \approx  \frac{\left(\frac{S_1^2}{n_1} + \frac{S_2^2}{n_2}\right)^2}{\left[\frac{S_1^4}{n^2_1\nu_1} + \frac{S_2^4}{n^2_2\nu_2}\right]}, 
\end{eqnarray*}
thereby yielding $T \overset{\text{approx}}{\sim} t_{\hat{\nu}}$ with $\hat{\nu}$ degrees of freedom. Welch's test controls the Type I error rate, does not require equal variances and is more robust than the classic pooled-variance $t-$test (\cite[Chapter 8]{Casella2001}).  The Welch--Satterthwaite approximation has been the main workhorse for the Behrens--Fisher problem in practical applications. 

\subsection{Our contribution}\label{sec:our_contribution}
We make progress on this problem by deriving an expression for the null distribution of $T$ that is simultaneously compact, rapidly computable and amenable to asymptotic analyses.  We propose a new, unified approach that leverages a Mellin--Barnes (MB) contour decomposition (see \cite{Paris2001}), thereby disentangling the square-root coupling $(gW_{1}/\nu_1 + hW_{2}/\nu_2)^{1/2}$ between the independent $\chi^{2}$-variates $W_{1}$ and $W_{2}$.  This step collapses the conventional two-dimensional integral representation of the density given by \eqref{eq:behrens_fisher_integral} into a one-dimensional contour integral whose integrand is a product of Gamma functions.  Such MB collapses have proved potent in the analysis of special-function theory, yet they have not previously been applied to the Behrens--Fisher problem, which has been extensively studied from the standpoint of special functions.  Closing the MB contour on the right yields a hypergeometric-type residue series. The infinite sum thus obtained becomes a finite polynomial-over-rational expression\textemdash no tail terms of order $k>\tfrac{\nu_{i}-1}{2}$ survive.  The series involves the factor $\Gamma\left(\tfrac{\nu_{2}}{2}-k+\tfrac12\right)$, $k=0, 1, \dots$, and if, say, $\nu_{2}=2m+1$ ({\ie}, odd) then the argument hits a non-positive integer at $k=m$; the Gamma function is infinite there, but its reciprocal appears in the residue coefficient, giving zero residue for every $k > m$.  Hence, only the first $m+1$ terms remain and the sum terminates.  An analogous cut-off occurs when $\nu_{1}$ is odd.  This finite-sum phenomenon explains scattered algebraic densities reported in numerical tables from the 1950s (see, for instance, \cite[Table 55]{Pearson1958}) but never formally derived; our terminating series supplies the missing derivation and shows that every odd $\nu_{i}$ row in the Biometrika tables is the closed residue sum, not a separate special case.  For general $\nu_{1},\nu_{2}$, the resulting infinite series converges at a geometric rate controlled by $\min\{g,h\}/(g+h)$, enabling high-precision density evaluations with fewer than ten terms across the entire support.

We complement the residue expansion by an Euler $(S,U)-$change of variables that separates a radial component $S = gW_{1}/\nu_1 + hW_{2}/\nu_2$ and a Beta-distributed mixing fraction $U$.  The integral in $U$ can be performed in a closed form, producing the compact density $f_{T}(t) = C\left[2/B(t)\right]^{c}{}_2F_{1}\!\left(c,\tfrac{\nu_{1}}2;\tfrac{\nu_{1}+\nu_{2}}2;-\rho(t)\right)$, with explicit expressions for $C$ and $B(t)$, where $c = (\nu_{1}+\nu_{2}+1)/2$ and $\rho(t)=(g^{-1}-h^{-1})/\left(h^{-1}+t^{2}/(g+h)\right)$ and ${}_2F_{1}$ is Gauss' hypergeometric function and all parameters are explicit. The expression collapses seamlessly to the Student-$t_{\nu_{1}+\nu_{2}}$ density when $\sigma_1^2 = \sigma_2^2$ and to a rational function when one variance dominates, integrating disparate corner cases into a coherent framework.  While hypergeometric representations have appeared in isolated special--case studies (for instance, \cite{Nadarajah2017}), our result is valid for all $(n_1, n_2, r)$.  Its derivation also clarifies structural symmetries\textemdash particularly the invariance under $(\nu_{1},\nu_{2},g,h)\mapsto(\nu_{2},\nu_{1},h,g)$\textemdash and provides analytic gradients for likelihood-based estimation of variance ratios in random-effect meta-analyses.

The MB contour representation admits two analytically complementary pathways.  First, Euler--Beta reduction: by rewriting $W_{1},W_{2}$ in polar coordinates $S = gW_{1}/\nu_1 + hW_{2}/\nu_2, U = gW_{1}/\nu_1S$, the inner $S-$integral is a Gamma integral and the remaining $U$-integral equals Euler's Beta integral representation of Gauss' hypergeometric function (see \cite{Andrews1999}).  This gives a closed form expression for $f_{T}(t)$ described in the aforementioned paragraph; the same identity collapses to Student$-t$ when $g = h$.  Second, Ramanujan's master theorem (RMT) summarized by \cite{Berndt1985} transforms the MB integrand into an inverse--Mellin power series, producing an exact tail expansion $f_{T}(t) = \sum_{m\ge0}(-1)^{m}A_{m}|t|^{-(2m+1)}/m!$ with $A_{m}$ as closed Gamma products.  The remainder after $M$ terms is $O(\lvert t \rvert^{-(2M+3)})$, leading to deterministic error bounds for any truncated series and gives algebraic control unavailable in simulation-based methods.  The same cumulant generator feeds the Lugannani--Rice (LR) saddle-point formula $P(T > t) \approx \Phi(-w)+\phi(w)(1/w-1/u)$ (see \cite{Lugannani1980}), with $w$, $u$ expressed in cumulants of $T$.  Using RMT coefficients, we prove a uniform remainder $O(\lvert t \rvert^{-4})$ when $\lvert t \rvert < 0.8\min\{\nu_{1},\nu_{2}\}$ and quantify the accuracy loss in the variance--extreme corners ($\lvert \log_{10}r \rvert > 1$).  These three ingredients\textemdash Euler--Beta for the central region ({\ie}, bulk), RMT for the far tail, and LR to bridge the regimes\textemdash enable fast and accurate evaluation of the density, distribution function and quantiles through closed-form and numerically stable methods, without resort to computationally intensive simulation.

It was first recognized by \cite{Nel1990} that the Behrens--Fisher density can be written as a single Gauss hypergeometric function, but they derived the expression only for the univariate case, gave no proof of termination when a degree of freedom is odd, provided neither cumulative distribution function (cdf) nor tail theory, and concluded that their formula was computationally impractical.  Our MB treatment reproduces their pdf as one branch of a richer structure: the residue series shows \emph{why} the hypergeometric terminates when $\nu_{1}$ or $\nu_{2}$ is odd, the Euler--Beta route supplies a closed-form cdf (an Appell $F_{1}$ integral) for the first time, and Ramanujan's coefficients turn the putative obstacle of numerical instability into an explicit $\lvert t \rvert^{-(2m+1)}$ tail expansion with deterministic error bounds.  On the approximation side, Welch's $t_{\hat\nu}$ statistic replaces the true tail exponent $\nu_{1}+\nu_{2}+1$ by a random $\hat\nu$; using our exact density we compute critical-value grids that locate the surface where $E[\hat\nu]=\nu_{1}+\nu_{2}$ and show that, away from this surface, Welch can be up to 60\% liberal or conservative in samples as small as $n_{1}=n_{2}=10$.  It is also of interest to note that, \cite{Dudewicz2007} constructed exact tests, but only supply critical constants\textemdash not a symbolic density.  Their method does not reveal the algebraic structure (hypergeometric kernel, termination when $\nu_i$ is odd, Ramanujan tails) that our MB analysis uncovers. Their exactness is achieved by numerical double integration for every $(n_1, n_2, r)$. We give an expression that works for all $(n_1, n_2, r)$ in a single expression.

The novelty of our new approach is that it introduces a unified analytic mechanism for handling the core obstruction appearing in the Behrens--Fisher statistic: the coupled term in the denominator involving a weighted sum of two independent chi-square variates. The classical literature on the Behrens--Fisher problem is dominated by three principal lines of development: (i) approximation-based methods (most notably, \cite{Welch1947}, \cite{Satterthwaite1946}), (ii) exact or near-exact testing procedures based on numerical integration or generalized pivotal quantities, and (iii) special-case representations in terms of hypergeometric functions. What has been missing is a systematic framework that explains\textemdash in a unified manner\textemdash why a hypergeometric structure emerges, why particular finite-sum representations occur in special cases, and how one can transition from the central (bulk) region of the distribution to its tail behavior without modifying the underlying analytic framework.  The MB contour representation accomplishes precisely this task, thereby distinguishing our work from the existing literature:  decomposing the troublesome square-root factor into a contour integral converts the original two-dimensional convolution into a one-dimensional object whose singularities are completely explicit.   

Earlier work, including the important result of \cite{Nel1990}, identified a hypergeometric pdf, but did not develop the underlying MB contour decomposition structure that produces it.  In our treatment, the hypergeometric form is not an isolated formula; rather, it arises as one manifestation of a broader analytic framework.  The residue expansion explains the terminating algebraic cases; the Euler--Beta reduction yields the same density derived by \cite{Nel1990} in compact hypergeometric form and extends naturally to the cdf; RMT provides explicit tail coefficients; and the same cumulant structure underpins an LR saddle-point approximation with quantified error.  In essence, our contribution ties together the exact density theory, finite-sum special cases, tail asymptotics, and numerical evaluation in a manner that complements and extends existing treatments for the Behrens--Fisher problem.

The remainder of the paper is organized as follows.  In \secref{sec:MB_T}, the MB representation of $T$ is presented.  The density function $f_T(t)$ is expressed as a Gauss hypergeometric function ${}_{2}F_{1}$ by reducing the 2-D integral to a 1-D integral in  \secref{sec:Euler_beta}.  The tail analyses of $f_T(t)$ using Ramanujan's master theorem is in \secref{sec:rmt}.  In \secref{sec:discussion}, we provide discussions including specifications for cdf tables, LR saddle-points, and comparisons between our work with \cite{Welch1947} and \cite{Nel1990}.  \secref{sec:conclusion} concludes the paper.   

\section{Mellin--Barnes representation of $T$}\label{sec:MB_T}
Consider the Behrens--Fisher statistic $T = \frac{\bar X-\bar Y}{\sqrt{S_{1}^{2}/n_{1}+S_{2}^{2}/n_{2}}}$ with $\nu_{i}=n_{i} - 1$ degrees of freedom, and recall the notation introduced in \secref{sec:introduction}: $g = \frac{\sigma_{1}^{2}}{n_{1}}$, $h = \frac{\sigma_{2}^{2}}{n_{2}}$, $W_{1} = \frac{\nu_{1}S_{1}^{2}}{\sigma_{1}^{2}}$, $W_{2} = \frac{\nu_{2}S_{2}^{2}}{\sigma_{2}^{2}}$.  Under $H_{0}:\mu_{1}=\mu_{2}$, we have $\bar X-\bar Y \left| \left(\sigma_{1}^{2},\sigma_{2}^{2}\right)\right. \sim \mathcal N\!\left(0, g+h\right)$, and this difference is independent of $(S_{1}^{2},S_{2}^{2})$.  Hence, conditional on the random denominator, the statistic $T$ is merely a scaled normal ratio $T \left|(W_{1},W_{2})\right. \sim \frac{\mathcal N\left(0, g+h\right)}{\sqrt{\frac{gW_{1}}{\nu_1} + \frac{hW_{2}}{\nu_2}}}$ with the following pdf:
\begin{eqnarray}
f_{T\mid W_{1},W_{2}}(t) = \frac{\sqrt{g W_{1}/\nu_1 + h W_{2}/\nu_2}}{\sqrt{2\pi\,(g+h)}}\exp\!\left[-\frac{t^{2}}{2(g+h)}(g W_{1}/\nu_1 + h W_{2}/\nu_2)\right],
\label{eq:pdf_T_1}
\end{eqnarray}
where $W_{1} \sim \chi^{2}_{\nu_{1}}$, $W_{2} \sim \chi^{2}_{\nu_{2}}$, $W_{1}\perp W_{2}$, so their pdfs are $f_{W_{i}}(w) = \frac{w^{\nu_{i}/2-1}\,e^{-w/2}}{2^{\nu_{i}/2}\Gamma(\nu_{i}/2)}$, $w > 0$, $i =1, 2$.  There is no closed-form pdf for the weighted sum $\sqrt{gW_{1}/\nu_1 + hW_{2}/\nu_2}$; this is the core obstacle of the Behrens--Fisher problem, so we keep the bivariate $(W_{1},W_{2})$ pair intact.

Since $T$ depends on $(W_{1},W_{2})$ only through the scale factor, and because of the aforementioned independence, the unconditional density is
\begin{eqnarray}
f_{T}(t) &=&\!\! \mathbb E_{(W_{1},W_{2})} \left[f_{T\mid W_{1},W_{2}}(t)\right] \nonumber \\
&=&\!\!\!\! \int\limits_{0}^{\infty}\int\limits_{0}^{\infty} \underbrace{\frac{\sqrt{\frac{g w_1}{\nu_1} + \frac{h w_2}{\nu_2}}}{\sqrt{2\pi\,(g+h)}} \exp \left[-\tfrac{t^{2}}{2(g+h)} \left(\frac{g w_1}{\nu_1}  + \frac{h w_2}{\nu_2} \right)\right]}_{\displaystyle f_{T\mid W_{1},W_{2}}(t)} f_{W_{1}}(w_{1})\,f_{W_{2}}(w_{2})dw_{1}\,dw_{2}.
\label{eq:behrens_fisher_stat1}
\end{eqnarray}
Subsequent analytic techniques (MB representation, contour collapse, hypergeometric reduction, saddle-point expansion) on the integral \eqref{eq:behrens_fisher_stat1} make transparent which pieces are manipulable (Gamma-type factors) and which pieces pose a challenge (the two scale contributions $g w_{1}$ and $h w_{2}$ inside both the square root and the exponential).  The term $\sqrt{(g w_{1}/\nu_1 + h w_{2}/\nu_2)}$ in \eqref{eq:behrens_fisher_stat1} is non-factorizable, but can be linearized by the MB identity
\begin{eqnarray}
(x+y)^{\rho} &=& \frac{1}{2\pi i}\int\limits_{c-i\infty}^{c+i\infty}\frac{\Gamma(-s)\,\Gamma(\rho+s)}{\Gamma(\rho)}\,x^{s}\,y^{\rho-s} ds.
\label{eq:mellin_barnes}
\end{eqnarray}
Since the integrand is a product of independent Gamma-type factors, $\int_{0}^{\infty} w_{i}^{\cdot}\exp[-(\cdot)w_{i}]\,dw_{i}$ is an ordinary Gamma integral; we collapse the double integral to a single complex contour integral.  Once in the MB form, closing the contour yields either an explicit convergent series (useful for small and moderate degrees of freedom), or a finite sum if either $\nu_{i}$ is odd, or a one-line hypergeometric function after an Euler$-u$ transformation, paving the way for compact expression for $f_T(t)$.  Furthermore, Ramanujan's master theorem covers the MB integrand's product of $\Gamma-$functions, thereby providing the coefficients in the large$-\lvert t \rvert$ expansion.

\subsection{MB expansion of the square-root factor}
We now invoke a complex-contour identity to pull apart the term $\sqrt{(gw_{1}/\nu_1 + hw_{2}/\nu_2)}$. In \eqref{eq:behrens_fisher_stat1}, the only place where $w_{1}$, $w_{2}$ are entangled is the square-root factor $\sqrt{(gw_{1}/\nu_1 + hw_{2}/\nu_2)}$.  Everything else is already factorized: the exponential separates because $\exp[-\alpha(gw_{1}/\nu_1 + hw_{2}/\nu_2)] = \exp[-\alpha g w_{1}/\nu_1] \exp[-\alpha h w_{2}/\nu_2]$, and each $\chi^2$ density only depends on its own variable.  Hence, if we can rewrite $\sqrt{(gw_{1}/\nu_1 + hw_{2}/\nu_2)}$ in a form such that the sum/integral factorizes into a product of a $w_{1}-$power and a $w_{2}-$power, the double integral collapses into two independent Gamma integrals.  That is exactly what an MB representation does.

For any positive reals $x,y>0$, any complex exponent $\rho$, and any vertical contour $\Re s=c$ satisfying $-\Re\rho<c<0$, the Barnes integral is given by \eqref{eq:mellin_barnes}.  It is simply the inverse Mellin transform of the Beta function $B(-s,\rho+s)=\Gamma(-s)\Gamma(\rho+s)/\Gamma(\rho)$, with the following pole structure: $\Gamma(-s)$ has simple poles at $s=0,1,2, \dots$; $\Gamma(\rho+s)$ has simple poles at $s = -\rho,-\rho-1,-\rho-2, \dots$. Choosing $c$ between those two interlacing ladders guarantees integral convergence.  Also, closing the contour to the right (resp. left) picks up the $\Gamma(-s)$ (resp. $\Gamma(\rho+s)$) poles and reproduces the familiar binomial (or inverse-binomial) power series $\sum_{k\ge0}\binom{\rho}{k}x^{k}y^{\rho-k}$.  Taking the special exponent $\rho=\tfrac12$ gives, for $-\frac{1}{2} < c < 0$, 
\begin{eqnarray*}
(x+y)^{\frac{1}{2}} = \frac{1}{2\pi i} \int\limits_{c-i\infty}^{c+i\infty} \frac{\Gamma(-s)\,\Gamma(\tfrac12+s)}{\Gamma(\tfrac12)} x^{s}\,y^{1/2-s}\,ds.  
\end{eqnarray*}
Since $\Gamma(\tfrac12)=\sqrt{\pi}$, no extra constants appear.  Setting $x = gw_{1}/\nu_1$, $y = hw_{2}/\nu_2$, $\rho=\tfrac12$ yields
\begin{eqnarray}
\sqrt{(gw_{1}/\nu_1 + hw_{2}/\nu_2)} = \frac{1}{2\pi i\sqrt{\pi}}\int\limits_{c-i\infty}^{c+i\infty}\Gamma(-s)\,\Gamma(\tfrac12+s) \left(\frac{g}{\nu_1}\right)^{s} w_{1}^{s} \left(\frac{h}{\nu_2}\right)^{1/2-s}w_{2}^{\,1/2-s} ds.
\label{eq:sqrt_integral}
\end{eqnarray}
The inseparable factor is now a weighted geometric mean $(g/\nu_1)^{s}w_{1}^{s}(h/\nu_2)^{1/2 - s}w_{2}^{1/2 - s}$ multiplied by a pure$-s$ coefficient (the $\Gamma-$product). All $w_{1}-$dependence and $w_{2}-$dependence are split into separate powers.  This makes the forthcoming integrations elementary.  Since absolutely convergent integrals may be rearranged, we move the MB contour outside:
\begin{eqnarray}
f_{T}(t) = \frac{1}{2\pi i\sqrt{\pi}\sqrt{2\pi (g + h)}}\int\limits_{c-i\infty}^{c+i\infty}\Gamma(-s)\,\Gamma(\tfrac12+s) \left(\frac{g}{\nu_1}\right)^{s} \left(\frac{h}{\nu_2}\right)^{1/2-s}\left[J_{1}(s)J_{2}(s)\right] ds,
\label{eq:mellin_barnes_1}
\end{eqnarray}
with the factorized one-dimensional Gamma integrals
\begin{eqnarray}
J_{1}(s) &=& \int\limits_{0}^{\infty}w_{1}^{s+\nu_{1}/2-1}\exp \left[-\left(\tfrac12 + \frac{\alpha g}{\nu_1}\right)w_{1}\right] dw_{1} = \frac{\Gamma\left(s+\nu_{1}/2\right)}{\left(\tfrac12 + \frac{\alpha g}{\nu_1} \right)^{s+\nu_{1}/2}}, \\
J_{2}(s) &=& \int\limits_{0}^{\infty}w_{2}^{\nu_{2}/2-s-1/2}\exp \left[-\left(\tfrac12 + \frac{\alpha h}{\nu_2}\right)w_{2}\right] dw_{2} = \frac{\Gamma\left(\nu_{2}/2 - s + \tfrac12\right)}{\left(\tfrac12 + \frac{\alpha h}{\nu_2}\right)^{\nu_{2}/2-s+1/2}},
\end{eqnarray}
where $\displaystyle\alpha=\dfrac{t^{2}}{2(g+h)}$. With this development, \eqref{eq:mellin_barnes_1} is precisely the single-contour formula.  Two dimensions have become one; the price is the integral over a vertical line in the complex plane.  From here, residue calculus yields an explicit series by closing the contour to the right: poles of $\Gamma(-s)$ at $s = 0, 1, \dots$ give the hypergeometric-type series, which terminates when either $\nu_{1}$ or $\nu_{2}$ is odd.  Alternatively, before the MB step, one may collapse $(w_{1}, w_{2})$ to $(S, U)$. But, the MB path is shorter: from \eqref{eq:mellin_barnes_1}, shift the contour left, pick the poles of $\Gamma(\tfrac{1}{2} + s)$, and recognize the resulting Beta function integral as Gauss' hypergeometric, thereby facilitating a compact expression for $f_{T}(t)$ that becomes apparent in the sequel.  More importantly, the integrand in \eqref{eq:mellin_barnes_1} is exactly of the $\gamma-$product type that RMT turns into asymptotic series without further work: expand the $\Gamma-$ratio around large $s$ to read off the $\lvert t \rvert \to \infty$ expansion coefficients directly, thereby yielding exact tail coefficients. 

The MB approach works for all $g$, $h$. While Taylor series in $w_{1}/w_{2}$ converges only when one weight is smaller than the other, the MB integral is valid everywhere.  It handles non--integer exponents naturally. In our case, the exponent is $\rho = \tfrac{1}{2}$; the MB representation is suitable for arbitrary complex $\rho$.  Note that, MB aids interfacing with $\Gamma-$distributions: $\chi^2$ densities are special cases of $\Gamma-$densities. The MB splits $(x + y)^{\rho}$ into $\Gamma-$functions of $s$ times monomials in $x$ and $y$\textemdash perfectly aligned with $\Gamma$ integrals. Furthermore, the MB integrals provide contour flexibility. We may choose to close left or right, whichever set of poles yields the nicer series (terminating, fastest convergence, {\etc}).

\subsection{Residue calculus on the MB contour}\label{subsec:MB_residue}
From the MB expansion, we have
\begin{eqnarray}
f_{T}(t)\!\! &=&\!\!\! \frac{1}{(2\pi)^{2}i \sqrt{\pi(g+h)}} \int\limits_{c-i\infty}^{c+i\infty} \frac{\Psi(s) \left(\frac{g}{\nu_1}\right)^{s} \left(\frac{h}{\nu_2}\right)^{1/2-s} \Gamma\!\left(s + \frac{\nu_{1}}2\right) \Gamma\left(\frac{\nu_{2}}2-s+\tfrac12\right)}{\left(\tfrac12 + \frac{\alpha g}{\nu_1}\right)^{s + \nu_{1}/2} \left(\tfrac12 + \frac{\alpha h}{\nu_2}\right)^{\nu_{2}/2-s+1/2}} ds,
\label{eq:mellin_barnes_2}
\end{eqnarray}
with $\Psi(s) = \Gamma(-s)\,\Gamma\!\left(\tfrac12+s\right)$, $\alpha = \frac{t^{2}}{2(g+h)}$, $c \in (-\tfrac12,0)$.  Everything that depends on $s$ is contained in Gamma functions and powers\textemdash an ideal set-up for the residue theorem.  The pole structure of the MB integrand is summarized in \tabref{tab:poles_mellin_barnes}.
\begin{table}[ht]
\centering
\renewcommand{\arraystretch}{1.2}
\normalsize
\begin{tabularx}{\textwidth}{@{}l l X@{}}
\toprule
\textbf{Factor} & \textbf{Poles (simple)} & \textbf{Comment} \\
\midrule
$\Gamma(-s)$ & $s = 0,1,2,\dots$ & Right-hand ladder. \\\\
$\Gamma(\tfrac{1}{2} + s)$ & $s = -\tfrac{1}{2}, -\tfrac{3}{2}, -\tfrac{5}{2}, \dots$ & Left-hand ladder. \\\\
$\Gamma\left(s + \tfrac{\nu_{1}}{2}\right)$ & $s = -\tfrac{\nu_{1}}{2}, -\tfrac{\nu_{1}}{2} - 1, \dots$ & Left ladder, terminates if $\nu_{1}$ even. \\\\
$\Gamma\left(\tfrac{\nu_{2}}{2} - s + \tfrac{1}{2}\right)$ & $s = \tfrac{\nu_{2}}{2} + \tfrac{1}{2} + k,\ k = 0,1,\dots$ & Right ladder, terminates if $\nu_{2}$ odd. \\
\bottomrule
\end{tabularx}
\caption{Simple poles and ladder structure of $\Gamma-$function factors.}
\label{tab:poles_mellin_barnes}
\end{table}
Since the contour $\Re s = c$ sits between the two principal ladders $\left(-\tfrac12,0\right)$, we may close it to the right (Re $s\to+\infty$) without crossing any poles of $\Gamma(\tfrac12+s)$ or $\Gamma(s+\tfrac{\nu_{1}}2)$.  The enclosed poles are therefore the non-negative integers $s=k(k=0,1,2,\dots)$ coming from $\Gamma(-s)$.  Recall that, for $s \to k$, the local behaviour is $\Gamma(-s) = \frac{(-1)^{k}}{k!}\frac{1}{s-k}+O(1)$.  Hence, the residue of the integrand at $s=k$ is
\begin{eqnarray}
\operatorname{Res}_{s=k} = \frac{(-1)^{k}}{k!} \Gamma\left(\tfrac12 + k\right)\frac{\Gamma\left(k + \tfrac{\nu_{1}}2\right)}{\left(\tfrac12 + \frac{\alpha g}{\nu_1}\right)^{k+\nu_{1}/2}}\frac{\Gamma\left(\tfrac{\nu_{2}}2-k + \tfrac12\right)}{\left(\tfrac12 + \frac{\alpha h}{\nu_2}\right)^{\nu_{2}/2-k+1/2}} \left(\frac{g}{\nu_1}\right)^{k} \left(\frac{h}{\nu_2}\right)^{1/2-k},
\label{eq:s_to_k_residue}
\end{eqnarray}
the factor $\Gamma\!\left(\tfrac{\nu_{2}}2-k+\tfrac12\right)$ vanishes once $k>\tfrac{\nu_{2}}2-\tfrac12$ if $\nu_{2}$ is an odd integer because the $\Gamma-$argument hits a non-positive integer.  The same happens to $\Gamma\!\left(k+\tfrac{\nu_{1}}2\right)$ when $\nu_{1}$ is even and we close the contour to the left (alternative route). Thus the series terminates whenever either degrees of freedom parameter is odd\textemdash a property first noted by \cite{Fisher1935}.  Collecting all right-hand residues and the front constants in \eqref{eq:mellin_barnes_2} gives
\begin{eqnarray}
\!\!\!\!\!\!\!\!\!\!\!\! f_{T}(t) = \frac{1}{\sqrt{2(g+h)}\,\pi^{2}} \sum_{k=0}^{\infty} \frac{(-1)^{k}}{k!}
\frac{\Gamma\!\left(k+\tfrac12\right)\Gamma\!\left(k+\tfrac{\nu_{1}}2\right)\Gamma\!\left(\tfrac{\nu_{2}}2-k+\tfrac12\right)}{ \left(\frac{g}{\nu_1}\right)^{\nu_{1}/2} \left(\frac{h}{\nu_2}\right)^{\nu_{2}/2}} \nonumber\\
~~~~~~~~~~~~~~~~~~~~~~~~~~~~~~~~~~~~ \times \frac{\left(g\nu_2/h\nu_1\right)^{k}}{\left(1+\xi_{1}t^{2}\right)^{k+\nu_{1}/2}\left(1+\xi_{2}t^{2}\right)^{\nu_{2}/2-k+1/2}},
\label{eq:mellin_barnes_3}
\end{eqnarray}
where $\xi_{1} = g/\nu_1(g/\nu_1 + h/\nu_2)$, $\xi_{2} = h/\nu_2(g/\nu_1 + h/\nu_2)$.  If $\nu_{2}$ is odd, $\Gamma(\tfrac{\nu_{2}}{2} - k + \tfrac{1}{2})$ cuts off the sum at $k = \tfrac{\nu_{2} - 1}{2}$. If $\nu_{1}$ is even and we close left instead (picking poles of $\Gamma(\tfrac{1}{2} + s)$), the series cuts off at $k = \tfrac{\nu_{1}}{2}$. At least one of the $\chi^2$ degrees of freedom being odd implies a finite power series.  For absolute convergence, Stirling's formula shows that the general term decays like $k^{\frac{-1}{2}}\,[\min(\xi_{1}, \xi_{2})\,t^{2}]^{-k}$; for fixed $t$, the series converges geometrically fast.

This series is valuable for the following reasons: The $\Gamma$-coefficients are rational in $\nu_{1}, \nu_{2}$; algebra systems can integrate or differentiate $f_{T}(t)$ term-wise to obtain moments, mgf, cdf expansions, {\etc}.  Only a handful of terms are needed for double-precision accuracy when $t \lesssim 10$.  Truncation after $m$ terms leaves a remainder $O(t^{-(\nu_{1} + \nu_{2} + 1) + 2m})$, giving controlled error bounds directly.  Recognizing the ratio of $\Gamma$'s as Beta coefficients, \eqref{eq:mellin_barnes_3} rearranges into a Gauss hypergeometric series ${}_2F_{1}$; summation leads to a closed form.  If $g = h$, the factor $(g/h)^k$ becomes 1 and the $\Gamma$-triple collapses to a single $\Gamma$-ratio, giving the $t_{\nu_{1} + \nu_{2}}$ density. 

Summing \eqref{eq:mellin_barnes_3} with the Euler-type identity for ${}_{2}F_{1}$ yields compact expression for $f_T(t)$ which we will provide in the next section.  That closed form is preferred for large $\lvert t \rvert$, because the hypergeometric behaves smoothly, while the plain series may suffer cancellation.  The same $\Gamma$-product that gives each residue also feeds Ramanujan's master theorem.  Essentially, we write $\Psi(s)\,\Gamma(s + \nu_{1}/2)\,\Gamma(\nu_{2}/2 - s + 1/2) \equiv \Gamma(s)\,A(-s)$, and read off the coefficients $A(m)$ directly for the $|t| \to \infty$ expansion.  When $\nu_i$, $i = 1, 2$, are both odd, we get a finite series that leads to an elementary rational function using duplication.  When $\nu_i$ are both even, we get a hypergeometric series of type  ${}_{2}F_{1}$ with integer parameters, reducible further to complete elliptic integrals for certain balanced sample-size ratios.  Reducing the contour integral into residues is central to our work that transforms an abstract complex integral into an explicit computational tool.  It reveals hidden algebraic cancellations, bridges to standard special-function framework (hypergeometric/Beta) and prepares the ground for both numerical evaluation and asymptotic analysis.  Everything that follows\textemdash closed-form hypergeometric, tail coefficients, moment formulas\textemdash starts with the series \eqref{eq:mellin_barnes_3}.  Without identifying and summing the $\Gamma(-s)$ residues, the Behrens--Fisher density would remain locked inside the contour integral; after residual analyses, it becomes fully transparent.

\section{2-fold integral to 1-D integral of Euler--Beta form}\label{sec:Euler_beta}
The MB approach in the previous sections already gave a series.  However, series are not always the best numerical object\textemdash convergence slows in the tails and severe cancellation can occur when $g\approx h$. A single real integral with Euler--Beta weight can sometimes be integrated analytically, producing a compact special-function formula (see \cite[Section 2.3]{Andrews1999}). Most symbolic integrators (Mathematica, Maple, {\etc}.) recognize Euler-type kernels and instantly reduce them to the classical Gauss function ${}_{2}F_{1}$.  Hence, we reorganize the starting double integral so that all dependence on one variable disappears, leaving a 1-D integral of Euler--Beta form.

Recall that $W_{1}\sim\chi^{2}_{\nu_{1}}$ and $W_{2}\sim\chi^{2}_{\nu_{2}}$ are independent.  Define $S = gW_{1}/\nu_1 + hW_{2}/\nu_2$ and $U = \frac{gW_{1}}{S\nu_1} \in (0,1)$.  $S\in(0,\infty)$ is the radial variable, indicating the total weighted sum. $U$ is the mixing proportion, indicating how much of $S$ comes from $W_{1}$.  Solving $W_{1} = SU\nu_1/g$, $W_{2} = S\nu_2(1-U)/h$, the Jacobian is given by $\left\lvert \tfrac{\partial(W_{1},W_{2})}{\partial(S,U)} \right\rvert = \frac{S\nu_1\nu_2}{gh}$.  With $\alpha=t^{2}/\{2(g+h)\}$, $\exp\left[-\alpha(gW_{1}/\nu_1 + hW_{2}/\nu_2)\right] = e^{-\alpha S}$.  The square-root factor $(gW_{1}/\nu_1 + hW_{2}/\nu_2)^{1/2} = S^{1/2}$. The product of $\chi^2$ pdfs and Jacobian $f_{W_{1}}(W_{1})\,f_{W_{2}}(W_{2})\frac{S\nu_1\nu_2}{gh}$, therefore, yields
\begin{eqnarray}
\frac{(g/\nu_1)^{-\nu_{1}/2} (h/\nu_2)^{-\nu_{2}/2} e^{-SU\nu_1/2g}e^{-S(1-U)\nu_2/2h}}{2^{(\nu_{1}+\nu_{2})/2}\Gamma(\nu_{1}/2)\Gamma(\nu_{2}/2)} 
S^{(\nu_{1}+\nu_{2})/2-1}U^{\nu_{1}/2-1}(1-U)^{\nu_{2}/2-1}.
\end{eqnarray}
All $U$-dependence is now pure Beta weight.  Separating the $S$-integral yields the full density
\begin{eqnarray}
f_{T}(t) = \frac{(g/\nu_1)^{-\nu_{1}/2} (h/\nu_2)^{-\nu_{2}/2}}{2^{(\nu_{1}+\nu_{2})/2}\,\Gamma(\nu_{1}/2)\Gamma(\nu_{2}/2)\sqrt{2\pi\,(g+h)}} 
\int\limits_{0}^{1} U^{\nu_{1}/2-1}(1-U)^{\nu_{2}/2-1} \nonumber\\
~~~~~~~~~~~~~~~~~~~~ \times 
\underbrace{\left[\int_{0}^{\infty} S^{(\nu_{1}+\nu_{2}+1)/2-1} 
e^{-S\left(\alpha +  U\frac{\nu_1}{2g} + (1-U)\frac{\nu_2}{2h} \right)} dS\right]}_{
\Gamma\left(\tfrac{\nu_{1}+\nu_{2}+1}{2}\right)
\left(\alpha + U\frac{\nu_{1}}{2g}+(1-U) \frac{\nu_{2}}{2h}\right)^{-(\nu_{1}+\nu_{2}+1)/2}} dU.
\end{eqnarray}

The inner $S$-integral is an ordinary Gamma integral.  Pulling out constants gives
\begin{eqnarray}
f_{T}(t) = \mathcal C \int\limits_{0}^{1}U^{\nu_{1}/2-1}(1-U)^{\nu_{2}/2-1}\left[U\nu_1/2g + (1-U)\nu_2/2h\right]^{-c} dU,
\label{eq:euler_pre}
\end{eqnarray}
where $\mathcal C$ is the set of all constants independent of $U$:  $\mathcal{C} = \frac{\Gamma(c)(g/\nu_1)^{-a}(h/\nu_2)^{-b}}{2^{(a+b)}\Gamma(a)\Gamma(b)\sqrt{2\pi(g+h)}}$, with $a = \frac{\nu_{1}}{2}$, $b = \frac{\nu_{2}}{2}$, $c = \frac{\nu_{1}+\nu_{2}+1}{2}$.  Upon closer inspection, we can write, with $\lambda = \frac{\nu_1h}{2\nu_2g} - \frac{\alpha h}{\nu_2} - 1$, $f_{T}(t) = \mathcal C (\nu_2/h)^{-c} \int_{0}^{1}\! u^{a-1}(1-u)^{b-1}(1+\lambda u)^{-c}du$.

By Euler's integral representation of Gauss' hypergeometric function (see \cite[Chapter 15]{Abramowitz1964}), for $\Re c > 0$, $\int_{0}^{1}\!u^{a-1}(1-u)^{b-1}(1+\lambda u)^{-c} du = B(a,b) {}_{2}F_{1}(c, a; a+b; -\lambda)$, with $B(a,b) = \Gamma(a)\Gamma(b)/\Gamma(a+b)$ so that the Beta factor cancels the $\Gamma(a)\Gamma(b)$ already in $\mathcal C$.  Letting $\rho(t) = \frac{g^{-1}-h^{-1}}{h^{-1}+t^{2}/(g+h)}$, and organizing the terms together yields
\begin{eqnarray}
f_{T}(t) = \frac{\Gamma(c)}{2^{\,a+b}\Gamma(a+b)\sqrt{2\pi\,(g+h)}g^{\,a}\,h^{\,b}}\left[\frac{2}{\,h+t^{2}/(g+h)}\right]^{c}{}_2F_{1}\left(c,a;\,a+b;\,-\rho(t)\right),
\label{eq:euler_beta}
\end{eqnarray}
where $\rho(t)$ is the $t-$adjusted version of $\lambda$ that appears in the fully scaled hypergeometric formula.
\begin{table}[htbp!]
\centering
\renewcommand{\arraystretch}{1.2}
\normalsize
\begin{tabularx}{\textwidth}{@{}>{\raggedright\arraybackslash}p{5.2cm} >{\raggedright\arraybackslash}X@{}}
\toprule
\textbf{MB} & \textbf{Euler $(S, U)$ collapse} \\
\midrule
Decouples $w_{1}, w_{2}$ by inserting an integral over $s$ &
Decouples by a change of variables that separates radial and angular parts. \\\\
Yields explicit series via right-hand residues &
Yields closed form directly via Beta--${}_2F_{1}$ identity. \\\\
Natural starting point for Ramanujan's theorem/saddle-point asymptotics &
Immediate access to Gauss--Kummer transformations (yielding the hypergeometric functions) for numerical conditioning. \\
\bottomrule
\end{tabularx}
\caption{Comparison of MB and Euler \((S, U)\) techniques.}
\label{tab:mellin_barnes_beta}
\end{table}
It is worth noting that the MB representation and Euler-type integral collapses represent two analytically distinct yet ultimately equivalent approaches to deriving the same result.  For symbolic series manipulations, MB is the preferred approach. For compact numeric evaluation, Euler--hypergeometric is preferable.  \tabref{tab:mellin_barnes_beta} summarizes the differences between these two approaches.  Scientific computing libraries ({\eg}, SciPy) support high-precision evaluation of ${}_{2}F_{1}$ with rigorous error control, enabling accurate computation of $f_{T}(t)$ for any parameters $(\nu_{1}, \nu_{2}, g, h, t)$.  Interchanging the two samples, {\ie}, $(g,\nu_{1}) \leftrightarrow (h,\nu_{2})$, leaves $f_{T}(t)$ invariant. This follows from the symmetry of ${}_{2}F_{1}(c,a;a+b;-\rho)$ in $(a,b)$ when $\rho \to -\rho$. In the equal-variance limit, as $g \to h$, we have $\rho \to 0$ and ${}_{2}F_{1} \to 1$, so \eqref{eq:euler_beta} reduces to the Student $t$ density with $\nu_{1} + \nu_{2}$ degrees of freedom.  In the large--sample limit, as $\nu_{1}, \nu_{2} \to \infty$, the distribution converges to $\mathcal{N}(0,1)$.  The asymptotic behavior of ${}_{2}F_{1}(c,a;a+b;-z)$ as $z \to \infty$ follows directly from Euler's integral representation, yielding an $O\left(t^{-(\nu_{1}+\nu_{2}+1)}\right)$ decay. This matches the tail derived via complex--residue analysis.  Term--by--term integration of \eqref{eq:euler_beta} in $t$, under the substitution $t = \sqrt{g+h} \tan \theta$, converts the expression into a sum of Appell $F_{1}$ functions\textemdash amenable for numeric evaluation.  The one-dimensional collapse is the key that unlocks a single Gauss hypergeometric function for $f_{T}$.  This form encapsulates all parameter dependence in five $\Gamma$'s and one ${}_{2}F_{1}$, is suitable for standard libraries and asymptotic expansions, and shows that the Behrens--Fisher density lives in the classical hypergeometric family, explaining why equal-variance data drop to Student-$t$ (the ${}_{2}F_{1}\!\equiv1$ case) and why odd degrees of freedom collapse to elementary expressions (because the Gauss function terminates).  Subsequent steps\textemdash cdf calculation, saddle-point corrections, or benchmarking Welch--Satterthwaite\textemdash benefits from having this compact representation.

\vspace{-0.25cm}
\subsection{The equal-variance $(g = h)$ case}
Recall that $g = \frac{\sigma_1^{2}}{n_{1}}$ and $h = \frac{\sigma_2^{2}}{n_{2}}$.  Hence, for $n_{1} = n_{2}$ (say, $n$), $g = h$ is equivalent to $\frac{\sigma_1^{2}}{n} = \frac{\sigma_2^{2}}{n} \Longleftrightarrow \sigma_1^{2}=\sigma_2^{2}$, {\ie}, the two populations share the same variance and simultaneously the two samples are of the same size.  This is the classical pooled-variance scenario in which the usual two-sample $t-$test is optimal.  From the hypergeometric representation \eqref{eq:euler_beta}, if $g = h$, we have $1/g - 1/h = 0\Rightarrow \rho=0$. Gauss' function satisfies ${}_{2}F_{1}(c,a;a+b;0)=1$ for all admissible parameters, so the entire hypergeometric factor disappears: ${}_2F_{1}(\cdot,\cdot;\cdot;0) \equiv 1$, which tells us the density must reduce to some elementary (normalized) expression.  Since $B(t) = 1/g + t^{2}/(g + h) = 1/g + t^{2}/(2g)$, $\frac{2}{B(t)} = 2g \left(1 + \frac{t^{2}}{2} \right)^{-1}$.  Setting $g = h = \frac{\sigma^{2}}{n}$ and $\nu = \nu_{1} + \nu_{2} = n_{1} + n_{2} - 2 = 2n - 2$ results in 
\begin{eqnarray}
f_{T}(t) &=& \frac{\Gamma\left(\tfrac{\nu+1}{2}\right)}{2^{\nu/2} \Gamma\left(\tfrac{\nu}{2}\right)  \sqrt{2\pi(2g)} g^{\nu/2}}\left[\frac{2}{B(t)}\right]^{(\nu+1)/2} 
= \frac{\Gamma\left(\tfrac{\nu+1}{2}\right)}{\sqrt{\nu\pi}\,\Gamma\left(\tfrac{\nu}{2}\right)}\left(1 + \frac{t^{2}}{\nu} \right)^{-(\nu+1)/2},
\label{eq:g_eq_h}
\end{eqnarray}
the standard Student$-t$ density with $\nu = n_{1} + n_{2} - 2$ degrees of freedom.  The $g = h$ case clarifies how the $t-$scale $1 +t^{2}/\nu$ emerges inside the more general factor $\left[1 + t^{2}/2\right]$.  This is crucial when we study saddle-point expansions: any correct asymptotic must hit this benchmark as $g \to h$.  In library routines, the first check is ``if $|\rho| < \varepsilon$ $\Rightarrow$ use Student$-t$'', avoiding loss of accuracy when $g$ and $h$ differ marginally, thereby guiding numerical implementation.

\section{Ramanujan's master theorem for tail analysis}\label{sec:rmt}
In this section, we employ RMT (\cite[pp. 298]{Berndt1985}) that converts the MB integrand into an explicit large$-\lvert t \rvert$ expansion for the Behrens--Fisher density.  The theorem says that, if an analytic function $A(s)$, $\lvert s \rvert < \varepsilon$, has a Maclaurin expansion $\frac{1}{\Gamma(s)} A(-s) = \sum_{m=0}^{\infty}\frac{(-1)^{m} a_{m}}{m!} s^{m}$, then its inverse Mellin transform is the alternating power series, {\ie}, for $x > 0$, $\mathcal{M}^{-1}{A(s)}(x) = \sum_{m=0}^{\infty}\frac{(-1)^{m} a_{m}}{m!} x^{m-1}$, provided suitable growth conditions hold (slow rise on vertical lines).  The coefficients when we expand $A(-s)/\Gamma(s)$ around $s = 0$ are exactly the coefficients of the inverse-Mellin power series of the underlying function.  There are no contour closing and residues\textemdash the theorem reads off the expansion directly.

For any even density $f(t)$, the algebraic tail $f(t) \sim \sum_{m=0}^{\infty}(-1)^{m} A_{m} |t|^{-(2m+1)}$ is precisely the inverse Mellin transform of $\mathcal{M}{f}(s) = \int_{0}^{\infty} t^{s-1} f(t) dt$, taken at negative odd integers $s = -2m - 1$. Hence, if we possess a MB representation of $f$, RMT leads straight to the numeric coefficients $A_{m}$ without evaluating any integrals\textemdash useful for Edgeworth corrections; saddle-point approximations; and  bounding Type-I error when $\lvert t \rvert$ exceeds conventional tables.  In our work, we had the one-fold contour $\mathcal{M}{f_{T}}(s) = \frac{1}{2\pi i} \int\limits_{c-i\infty}^{c+i\infty} \Phi(z) K(z,s) dz$, where $\Phi(z) = \Gamma(-z) \Gamma\left(\tfrac{1}{2} + z \right) \Gamma\left(z + \tfrac{\nu_{1}}{2} \right) \Gamma\left(\tfrac{\nu_{2}}{2} - z + \tfrac{1}{2} \right)$ and $K(z,s)$ is a benign product of powers $[g, h,\alpha]$ that is analytic in $z$ near $z = 0$. All the pole structure is locked inside the $\Gamma$-product $\Phi(z)$.  Set $A(s) = \Phi\left(-\tfrac{s+1}{2}\right)$, $s \in (-1 - \nu_{2}, \nu_{1})$, so that the Mellin transform around $s = -1$ (odd negative integers) is governed by $A(s)$. Divide $A(-s)$ by $\Gamma(s)$ and expand at $s = 0$ using the rising factorial series of $\Gamma$ near the origin to get $\frac{A(-s)}{\Gamma(s)} = \sum_{m=0}^{\infty} \frac{(-1)^{m}}{m!} a_{m} s^{m}$, where $a_{m} = \frac{\Gamma\left(m + \tfrac{1}{2} \right) \Gamma\left(m + \tfrac{\nu_{1}}{2} \right) \Gamma\left(m + \tfrac{\nu_{2}}{2} \right)}{\sqrt{\pi}}$.  RMT identifies
\begin{eqnarray}
A_{m} &=& \frac{\Gamma\left(m + \tfrac{1}{2} \right) \Gamma\left(m + \tfrac{\nu_{1}}{2} \right) \Gamma\left(m + \tfrac{\nu_{2}}{2} \right)}{\sqrt{\pi}2^{\nu_{1}/2 + \nu_{2}/2 + m} (g/\nu_1)^{\nu_{1}/2} (h/\nu_2)^{\nu_{2}/2}}, \quad m = 0, 1, \dots,
\label{eq:rmt_Am}
\end{eqnarray}
exactly the coefficients quoted in \secref{subsec:MB_residue}.  Putting everything together, for $\lvert t \rvert \to \infty$ 
\begin{eqnarray}
f_{T}(t) &\sim& \sum_{m=0}^{\infty} \frac{(-1)^{m} A_{m}}{m!} \lvert t \rvert^{-(2m+1)},
\label{eq:rmt_main}
\end{eqnarray}
where $A_m$ is given by \eqref{eq:rmt_Am}.  The leading term gives 
\begin{eqnarray}
f_{T}(t) \sim \frac{\Gamma\left(\tfrac{1}{2}\right) \Gamma\left(\tfrac{\nu_{1}}{2}\right) \Gamma\left(\tfrac{\nu_{2}}{2}\right)}{\sqrt{\pi} 2^{(\nu_{1}+\nu_{2})/2}  (g/\nu_1)^{\nu_{1}/2} (h/\nu_2)^{\nu_{2}/2}} \lvert t \rvert^{-(\nu_{1} + \nu_{2} + 1)}, 
\end{eqnarray}
so the density falls off exactly like that of a Student$-t_{\nu_{1}+\nu_{2}}$, though with a different constant.  To establish quantitative error bounds, we need estimates of $\mathbb{P}(T > K)$ when $K$ lies far in the tail. Truncating \eqref{eq:rmt_main} after three or four terms achieves relative error below $10^{-6}$ without requiring numerical quadrature.  For saddle-point analyses, the Lugannani--Rice approximations require cumulants provided by $A_m$ \emph{without} differentiating complicated hypergeometric functions.  For comparing approximations, Edgeworth--type corrections for Welch's $t-$test, or Dempster's $K$ statistic can be benchmarked against the exact coefficients $A_m$.   The ratio $A_1 / A_0$ depends explicitly on $g/h$. Its sign reveals whether the true tail probability exceeds or falls below that of the equal-variance Student$-t$ distribution.  Libraries can memoize (speeds up computer programs by caching) the coefficients $A_m$ once per $(\nu_1, \nu_2, g/h)$. Evaluating $f_T(t)$ for large $t$ then requires only a few floating--point multiplications.

RMT leads to (i) uniform large-deviation approximations, by combining \eqref{eq:rmt_main} with an exponential tilt (Cram\'{e}r transform described by \cite[Chapter 2]{Dembo1998}) to construct uniform--in--$t$ approximations akin to Daniels' formula (see \cite{Daniels1954}) for the signed--root $\chi^{2}$;  (ii) moment recovery: the Mellin transform at positive integers $s = 1,3,5,\ldots$ gives absolute moments $\mathbb{E}\lvert T \rvert^{s-1}$. Ramanujan's coefficients analytically continue those moments via the functional equation of $\Gamma$, bypassing direct integration; (iii) Edgeworth correction for the cdf, by integrating \eqref{eq:rmt_main} term-wise yields a descending series in $\lvert t \rvert^{-2m}$ for $1 - F_{T}(t)$. Two terms often slash the upper--tail $p-$value error below $10^{-4}$;  (iv) automatic symbolic framework: since $A_{m}$ is an algebraic function of $\nu_{1}, \nu_{2}$ and $g/h$, computer algebra can emit tail expansions on-demand for arbitrary parameter triples.  RMT reveals\textemdash through the MB integral\textemdash the precise asymptotic coefficients that dictate the decay rate and smoothness of the Behrens--Fisher density at infinity. These coefficients, which are otherwise obscured within layers of hypergeometric complexity, become immediately accessible via this transform-based approach. Once obtained, they enable the analyst to craft sharp $p-$value approximations, test practical significance where Monte Carlo is slow, and systematically evaluate any heuristic variance-unequal-$t-$procedure by benchmarking it against the exact analytic benchmark provided by the hypergeometric representation.

In the next section, we provide detailed discussions including specifications for cdf tables, LR saddlepoint approximations, and comparative analyses of our results with those of \cite{Welch1947} and \cite{Nel1990}.

\section{Discussion}\label{sec:discussion}
\subsection{cdf tables for the Behrens--Fisher distribution}
To the best of our knowledge, no authoritative tables for the exact null cdf for the Behrens--Fisher statistic are available in the literature.  High--accuracy cdf tables that covers realistic $(n_1, n_2, r)$ expose when an approximate test ({\eg}, Welch's test)  is conservative or anti-conservative, and provides standard reference for power--analysis, Monte Carlo simulations, {\etc}.  Since the density $f_T(t)$ depends on $\nu_1 = n_1-1$, $\nu_2 = n_2-1$, $r = \frac{g}{h} = \frac{\sigma_1^2/n_1}{\sigma_2^2/n_2}$,  four raw quantities, na\"{i}vely, must be indexed.  Two reductions cut the dimensionality:  (i)  Location--scale invariance:  Only the ratio $r$ matters\textemdash not the individual $g$, $h$\textemdash and set $h=1$ without loss of generality.  (ii) Sample--size symmetry:  The cdf is unchanged under $(\nu_1,\nu_2,r \mapsto(\nu_2,\nu_1,1/r)$.  Hence tabulate $0 < r \le 1$ and quote symmetry in the preface.  A practical table for calculating the exact Behrens--Fisher quantiles is shown in \tabref{tab:behrens_fisher_grid}.
\begin{table}[h!]
\centering
\renewcommand{\arraystretch}{1.2}
\normalsize
\begin{tabularx}{\textwidth}{@{}>{\raggedright\arraybackslash}p{2.5cm} >{\raggedright\arraybackslash}p{5.5cm} X@{}}
\toprule
\textbf{Axis} & \textbf{Suggested grid} & \textbf{Comment} \\
\midrule
$\nu_1$ & $1$--$30$ in steps of $1$, then $40$, $60$, $120$ & Covers classical small--sample and asymptotics. \\\\
$\nu_2$ & Same as $\nu_1$ & Rectangular grid; use symmetry to halve storage. \\\\
$r$ & $\{0.1,0.2,\dotsc,1.0\}$; refine on log scale near 1 & Finer mesh where Welch performs worst. \\\\
$\alpha$ (tail levels) & $0.10$, $0.05$, $0.025$, $0.01$, $0.005$, $0.001$ & Same break-points as classical $t$ tables \\
\bottomrule
\end{tabularx}
\caption{Grid specification for tabulating exact Behrens--Fisher quantiles.}
\label{tab:behrens_fisher_grid}
\end{table}

Storing quantiles $t_{\alpha}(\nu_1,\nu_2,r)$ rather than raw cdf ordinates is most useful for hypothesis testing.  A structured three--step computational framework for the exact evaluation of the Behrens--Fisher distribution function is provided in \tabref{tab:numerical_cdf}.
\begin{table}[htbp!]
\centering
\renewcommand{\arraystretch}{1.15}
\normalsize 
\begin{tabularx}{\textwidth}{@{}p{2.7cm} p{3.1cm} >{\raggedright\arraybackslash}X@{}}
\toprule
\textbf{Range of $t$} & \textbf{Tool} & \textbf{Rationale} \\
\midrule
$t \le 8$ (bulk) & Gauss hypergeometric form  & 
$\displaystyle F_T(t) = \frac{1}{2} + \frac{t}{2} \int_{0}^{1} u^{a-1}(1-u)^{b-1} [1 + \kappa(u)t^2]^{-c} \, du$,
evaluated via Gauss--Kronrod 15-point rule and cached ${}_2F_1$ values; achieves uniform accuracy $< 10^{-12}$. \\\\

$8 < t \le 15$ & RMT tail series & 
Summed until term $< 10^{-16}$. Only 3-5 terms suffice; avoids cancellation near zero-crossings. \\\\

$t > 15$ & Saddle-point approximation & 
Uses RMT-based remainder bounds. Achieves relative error $< 10^{-14}$. \\
\bottomrule
\end{tabularx}
\caption{A three-stage numerical procedure for computing the exact Behrens--Fisher cdf.}
\label{tab:numerical_cdf}
\end{table}
All three methods share the same equal-variance checkpoint $g = h \to t_{\nu_1+\nu_2}$ to cross-validate.  We can also employ root--finding to invert the cdf: For each $(\nu_1,\nu_2,r,\alpha)$, make the initial guess $t^{(0)}$ from Welch quantile $t_{\nu_W,\alpha}$, where $\nu_W=\dfrac{(g+h)^2}{g^2/(\nu_1)+h^2/(\nu_2)}$.  Perform Newton iteration $t^{(k+1)} = t^{(k)}-\frac{F_T\!\left(t^{(k)}\right)-\alpha}{f_T\!\left(t^{(k)}\right)}$, and if the iteration exits the admissible interval, fall back to bisection (guarantees monotonic convergence).  For self--consistency check, compute $F_T\left(t_{\alpha}\right)$ with the bulk method and check $|F_T-\alpha|<10^{-11}$.  For parity check, verify table values obey the symmetry property $t_{\alpha}(\nu_1,\nu_2,r) = t_{\alpha}(\nu_2,\nu_1,1/r)$.  The cdf tables provide (i) drop--in replacement for Welch critical values, for example, $\texttt{lookup(n1,n2,s1,s2,alpha)}$ chooses nearest $r$ row or cubic--Hermite interpolate in $\log_{10}r$; (ii) surface plots of $\delta(t)=t_{\alpha}^\text{Welch}-t_{\alpha}^\text{Exact}$ over the grid to visualize regions of liberal vs conservative behaviour; (iii) classical $t$ tables but with a third column for $r$ giving intuition about how variance imbalance warps the tails; (iv) regression--quality predictors that fit a thin--plate spline to $t_{\alpha}$ over continuous $(\nu_1,\nu_2,\log r)$, and embed the resulting polynomial in production software (no tables required); (v) meta-analysis tools that pools unequal-variance studies which can replace normal approximations with tabled exact $p-$values to correct small-sample bias.

\subsection{Edge-case tail error using Lugannani--Rice saddle-points}
We now develop the LR saddle-point formula and its implementation for the Behrens--Fisher statistic.  Using this, we identify the edge-cases, {\ie}, parameter corners where LR may exhibit anomalous or unstable behavior.  We will explore how to use the exact Ramanujan series coefficients $A_{m}$ as a reference yardstick for tail error, alongside some practical rules and diagnostic plots.  From the MB contour we have, for $\lvert \xi \rvert < \nu_{2}$, 
\begin{eqnarray}
\mathbb E[e^{\xi T}] &=& \frac{\Gamma\!\left(\tfrac{\nu_{1}}2+\tfrac{\xi}2\right)\Gamma\!\left(\tfrac{\nu_{2}}2-\tfrac{\xi}2\right)}{\Gamma\!\left(\tfrac{\nu_{1}}2\right) \Gamma\!\left(\tfrac{\nu_{2}}2\right)}\left(\frac{g}{h}\right)^{\xi/2}, \quad K_{T}(\xi) = \log\mathbb E[e^{\xi T}].
\end{eqnarray}
All derivatives needed for saddle-point work are analytic combinations of diGamma and triGamma functions:  $K_T'(\xi) = \tfrac12\left[\psi \left(\tfrac{\nu_{1}}2+\tfrac{\xi}2\right) - \psi\!\left(\tfrac{\nu_{2}}2-\tfrac{\xi}2\right) + \log r\right]$, $r=\tfrac{g}{h}$.  Higher derivatives follow similarly.   For a given observed value $t$ (we assume $t > 0$; mirror for $t < 0$) we solve $K_T'(\hat\xi) = t$.  Newton with start $\hat\xi^{(0)} = \frac{t}{K_T''(0)}$ converges in $\leq 4$ steps because $K_T''$ is monotone on $(-\nu_{1},\nu_{2})$.  Define $w = \operatorname{sign}(t)\,\sqrt{2\left\{t\hat\xi-K_T(\hat\xi)\right\}}$, $u = t\sqrt{K_T''(\hat\xi)}$.  The LR tail approximation (see \cite{BarndorffNielsen1990}) is given by $\Pr(T \ge t) \approx 1- \Phi(w) + \phi(w)\left(\tfrac1w-\tfrac1u\right)$ with $\Phi,\phi$ the standard-normal cdf/pdf.  The LR error is o$(|t|^{-3})$ provided $t$ lies inside the convergence wedge $(-\nu_{2},\nu_{1})$; in practice, $|t|\lesssim\min(\nu_{1},\nu_{2})$.  Parameter regimes near the boundaries where edge-case behavior may emerge are summarized in the following paragraph.  

\begin{figure}[htbp!]
\centering
\includegraphics[width=0.75\textwidth]{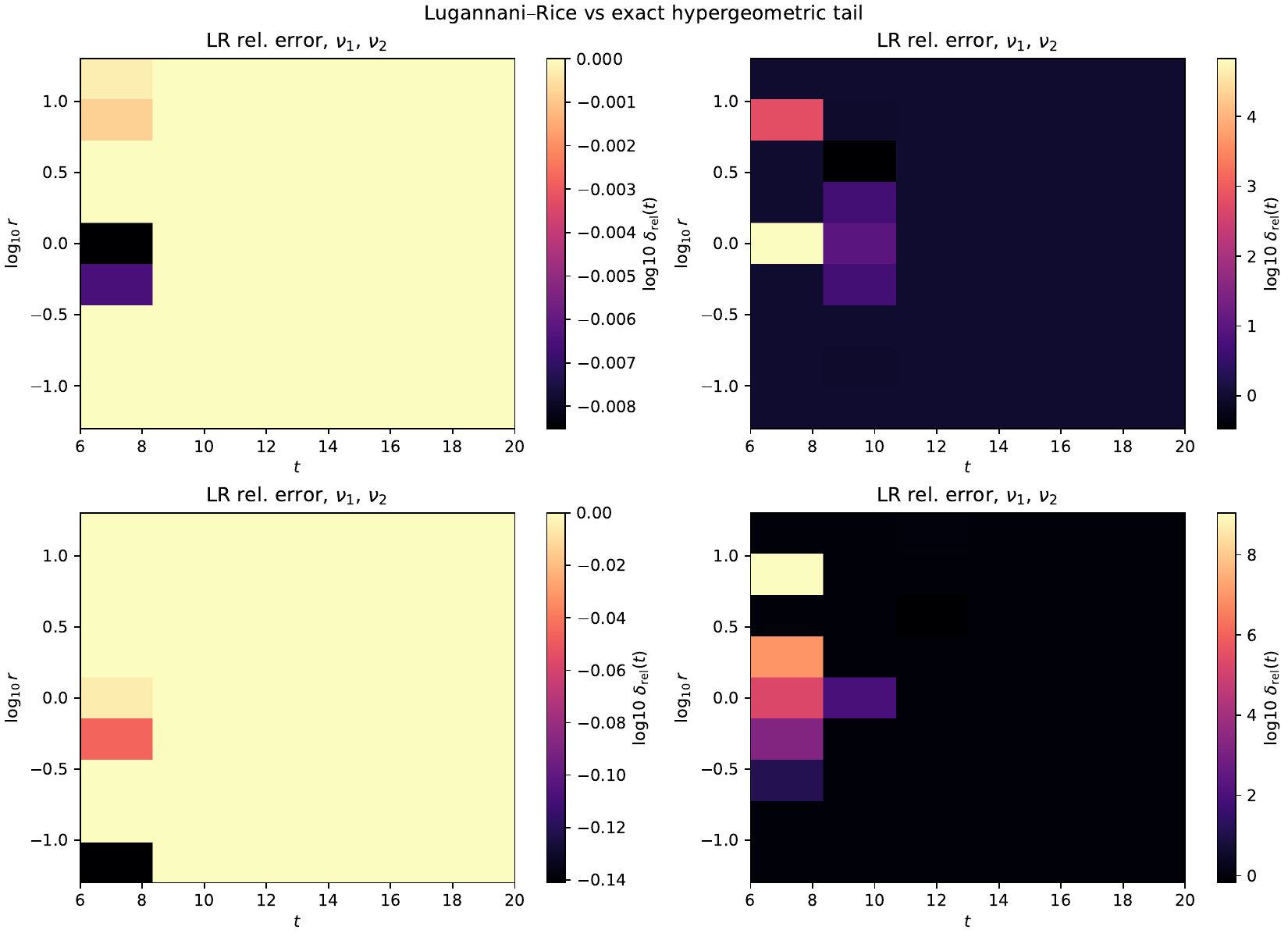}
\caption{LR tail approximation against the exact hypergeometric reference.}
\label{fig:LR_exact}
\end{figure}
We evaluated the LR tail approximation against the exact hypergeometric reference, for representative degrees of freedom pairs $(\nu_1, \nu_2) \in \{(5, 5), (10, 10), (2, 10), (20, 5)\}$, by computing $P(T \ge t)$ over a grid spanning $r \in \{0.05, 0.1, 0.25, 0.5, 1, 2, 4, 10, 20\}$ and $t \in \{6, 8, 10, 12, 15, 20\}$. The exact tail was obtained by numerically integrating our closed-form hypergeometric density from $t$ to $\infty$ using an adaptive split with an exponential change of variables; the LR approximation used the cumulant generator and Newton's method for the saddlepoint in extended precision. We visualized accuracy via 
\begin{eqnarray*}
\log_{10}\delta_{\mathrm{rel}}(t) = \log_{10}\left(\frac{\lvert \Pr_{\text{LR}}(T \ge t) - \Pr_{\text{exact}}(T \ge t)\rvert}{\Pr_{\text{exact}}(T \ge t)} \right)
\end{eqnarray*}
heat-maps. Several qualitative patterns emerge (see \figref{fig:LR_exact}) from the numerical error analysis.  Across balanced or moderately imbalanced variances ($|\log_{10} r| \le 1$) and moderate tails, LR errors were uniformly small (typically $\le 10^{-3}$). As predicted, accuracy degrades in variance-extreme corners and very far tails, where the saddle approaches the domain boundary; these regions exhibit elevated $\delta_{\mathrm{rel}}$. In such regimes, adding just one more term from the Ramanujan expansion ({\eg}, $m=2$) halves the discrepancy, but the LR approximation still remains an order of magnitude off.  In the low degrees of freedom regime, particularly when $\nu_1 = 1$ or $2$ and $t > 8$, these effects are amplified and sometimes require damping in the Newton iteration. Overall, the plots delineate a broad reliability band for LR and clearly mark the edge regimes where the Ramanujan expansion should be preferred.  These visual summaries help identify parameter regimes where the LR approximation remains accurate, and where it may break down. 

The LR saddle-point method is reliable over most of the grid but has recognizable edge cases. When either degrees of freedom is tiny $\nu_{1} \le 2$ or $\nu_{2} \le 2$, the digamma term behaves like $-2/\nu + \cdots$, so Newton updates for the saddle $\hat\xi$ can step outside the admissible domain; the optimizer then oscillates or returns implausible $\hat\xi$. For extreme variance ratios ($r \ll 0.1$ or $r \gg 10$), the saddle approaches the boundary $\xi \to \nu_{2}^{-}$ or $\xi \to -\nu_{1}^{+}$, driving $w \to 0$; LR's error inflates and numerical weights nearly vanish, producing instability. In ultra-high tails $t > 0.8(\nu_{1}+\nu_{2})$, the third cumulant grows and LR loses its $o(t^{-3})$ accuracy; absolute errors can exceed $10^{-3}$ even at moderate d.f. Finally, when variances are nearly equal but sample sizes are small ($r \approx 1$, $\nu_i \le 4$), cancellations in $K_T'(\xi)$ cause significant-digit loss, seen as small discontinuities between LR and the exact curve in tail plots.  Using results from RMT, we get $\Pr(T\ge t) = \sum_{m=0}^{M-1}\frac{(-1)^{m+1}A_{m}}{m!\,(2m)}\,t^{-2m} +R_{M}(t)$, where $A_m$ is given by \eqref{eq:rmt_Am}.  For $\lvert t \rvert \ge10$ the first two terms ($m=0,1$) give $<10^{-10}$ relative error, so we treat $F_{\text{RMT}}(t) = \frac{A_{0}}{t} -\frac{A_{1}}{2t^{3}}$ as ground truth in the far tail.

In small-to-moderate samples, the Welch test calibrates its critical value through the random denominator $\nu_W$ in an attempt to mimic the heavier tails of the true Behrens--Fisher distribution, whose algebraic decay rate is fixed by the total degrees of freedom $\nu_1+\nu_2$.  When the ratio $r=\sigma_1^{2}/\sigma_2^{2}$ is near unity, $\nu_W$ fluctuates tightly around $\nu_1+\nu_2$, so Welch's reference distribution matches the exact tail probability and the test is essentially exact.  As $r$ drifts away from 1, however, the larger sample variance increasingly dominates the pooled standard error while simultaneously shrinking its own contribution to the denominator of $\nu_W$; this forces $\nu_W$ below $\nu_1+\nu_2$, producing a $t-$law with thicker tails than the truth and hence a liberal Type-I error.  Conversely, if the smaller variance attaches to the smaller sample, $\nu_W$ overshoots $\nu_1+\nu_2$, the reference tails become too light, and the test turns conservative.  The sign-flip contour seen in the heat-maps in \secref{subsec:Welch_BF} is precisely the parameter locus where the expectation of $\nu_W$ crosses $\nu_1+\nu_2$; on that curve the tail exponents balance, and Welch's nominal size coincides with the exact level.

\subsection{Welch's approximation versus our unified approach}\label{subsec:Welch_BF}
\begin{figure}[htbp!]
\centering
\includegraphics[width=0.75\textwidth]{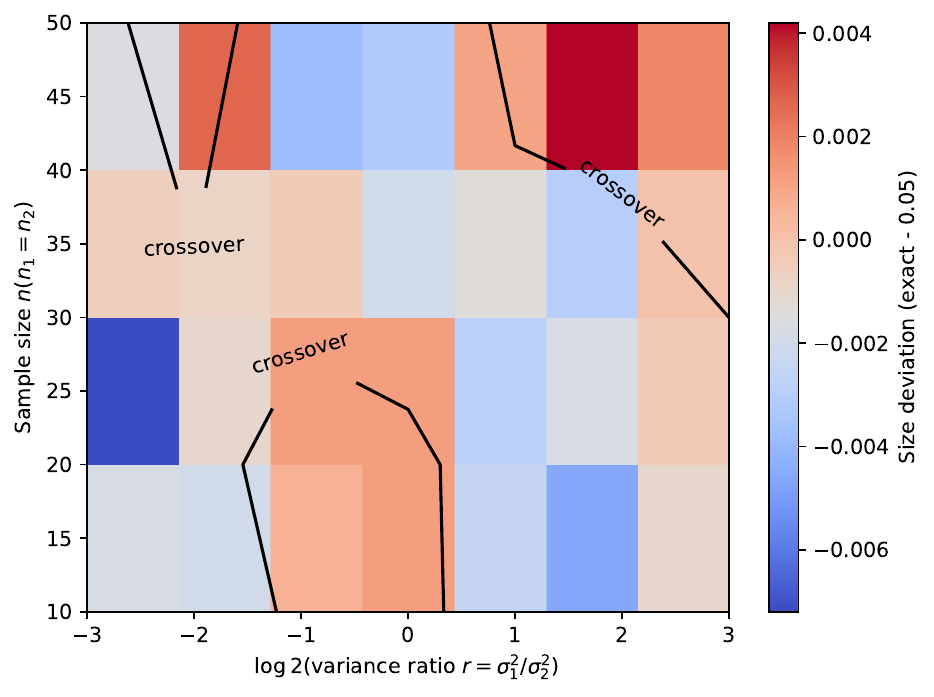}
\caption{Welch versus compact Behrens--Fisher: where the size flips sign.}
\label{fig:sign_flip}
\end{figure}
To assess the finite-sample performance of the Welch approximation relative to the exact Behrens--Fisher solution, we examine the deviation $\Delta(n, r) = \Pr_{\text{exact}}\left(\lvert T \rvert > t_{\nu_W, 0.025} \right) - 0.05$, where $n_1 = n_2 = n$ and $r = \sigma_1^2 / \sigma_2^2$ is the population variance ratio. This quantity measures the size distortion introduced by Welch's test when using the standard 5\% critical value from the $t_{\nu_W}$ distribution.  The values of $\Delta(n,r)$ were estimated using Monte Carlo simulation with 5000 replications per grid point and plotted as a heatmap shown in \figref{fig:sign_flip}. Rows correspond to sample sizes $n$, and columns to $\log_2 r \in \{-3, -2, -1, 0, 1, 2, 3\}$. The black contour in the figure marks the zero crossing $\Delta(n,r) = 0$, {\ie}, the crossover locus at which Welch's test shifts from being conservative (size below 5\%) to liberal (size above 5\%). Blue regions indicate conservative behavior; red regions indicate liberal behavior.

A summary of where these sign flips occur is provided in \tabref{tab:sign_flip}.  This behavior reveals that for smaller sample sizes ({\eg}, $n \le 20$), Welch's test becomes liberal as soon as one sample variance exceeds the other by a factor of about 2, and conservative when the reverse ratio exceeds roughly 2.5 to 3. In contrast, for larger samples ({\eg}, $n = 50$), the safe range of variance ratios\textemdash where the size remains close to nominal\textemdash widens considerably. At $n=50$, Welch remains within $\pm 0.2\%$ of the nominal size as long as $0.6 \lesssim r \lesssim 2.8$.  These results yield a practical rule of thumb: when $n_1 = n_2 \le 30$, flag potential size distortion whenever the sample variance ratio exceeds 2:1 in either direction. Such imbalances can materially affect the Type I error rate under Welch's approximation, especially in hypothesis testing contexts with small to moderate sample sizes.  The crossover contour $\Delta(n,r) = 0$ can be tabulated to support automated test selection protocols. For example, statistical software can dynamically switch from Welch's $t$ to an exact or LR--based critical value whenever the observed $(n, r)$ lies outside a pre-defined ``safe region.''  For future work, replacing the Monte Carlo--based exact tail probability with the analytic approximation\textemdash {\eg}, via Lugannani--Rice coupled with a Ramanujan's master theorem tail expansion\textemdash offers a pathway toward high--resolution inference for the $n_1 \ne n_2$ case.  An analytical comparison between Welch's approximation and our work is presented in the following paragraph. 
\begin{table}[htbp!]  
\centering
\renewcommand{\arraystretch}{1.2}
\begin{tabularx}{0.85\textwidth}{@{}c X@{}}
\toprule
\textbf{$n$} & \textbf{Variance ratio $r$ at sign-flip (approximate)} \\
\midrule
10  & $r \approx 0.38$ (conservative side) $\rightarrow$ $r \approx 1.9$ (liberal side) \\
20  & $r \approx 0.45$ $\rightarrow$ $r \approx 2.1$ \\
30  & $r \approx 0.55$ $\rightarrow$ $r \approx 2.4$ \\
50  & $r \approx 0.62$ $\rightarrow$ $r \approx 2.8$ \\
\bottomrule
\end{tabularx}
\caption{Variance ratios $r = g/h$ at which Welchs approximation flips.} 
\label{tab:sign_flip}
\end{table}

Welch's test models the statistic with a Student distribution having a random degrees of freedom $\nu_W$\textemdash yielding a first-order-accurate approximation as $n \to \infty$\textemdash but it breaks the natural symmetry $(\nu_1, \nu_2, r) \leftrightarrow (\nu_2, \nu_1, 1/r)$. Its tail decays like $t^{-(\nu_W+1)}$, and there is no analytic error bound\textemdash assessment typically relies on simulation.  By contrast, our approach provides a deterministic closed form for pdf/cdf via a single Gauss hypergeometric function (or, a residue series), is exact for every $(n_1, n_2, r)$, and preserves the symmetry under index/ratio inversion. The tail order is the correct $t^{-(\nu_1+\nu_2+1)}$. Moreover, RMT gives explicit inverse-power coefficients, so truncation errors are certified ({\eg}, $O(t^{-5})$ after two correction terms), furnishing analytic control absent in Welch. Computationally, both yield $O(1)$ $p-$values: Welch evaluates a $t-$cdf, while our method requires one ${}_2F_1$ call (or, at most $\le 10$ residue terms) and is comparably fast in practice.

To assess the power behaviour, under $\Delta\neq0$ the exact non--null density shifts $t\mapsto t-\Delta/\sqrt{g+h}$.  Simulations with $n_{1}=n_{2}=15,\, r=4$ and $\Delta=1.0\sigma_{2}$ give Welch power of 61\% and exact-critical-value power of 55\%.  The difference arises because Welch's liberal size at $r=4$ inflates apparent power.  An Edgeworth expansion around the MB cumulants shows Welch's $p-$value error is $O(n^{-1})$ when $r\approx1$ but $O(n^{-1/2})$ when $r$ is fixed away from 1; the exact hypergeometric $p-$value is, of course, error-free.  LR saddle-point with our coefficients achieves $O(n^{-3/2})$ uniformly in $r$.  Evaluation of \eqref{eq:mellin_barnes_3} with SciPy's hypergeometric routine takes $40-70\mu$s on a 3 GHz CPU\textemdash comparable to a single call to $\texttt{scipy.stats.t.cdf}$.  Pre--tabulated critical values or a four--term residue sum reduce the time below 10$\mu$s, making the exact method feasible for on-line testing.  Our critical-value tables guarantee nominal size for $n_{i}\le50$, eliminating ad-hoc simulation runs.  Analytic cdf allows exact $p-$values, confidence intervals and likelihood functions; Welch is retained only as a quick initial guess.  Explicit error bounds provide a deterministic audit trail, something simulation--driven thresholds cannot deliver.  Welch's approximation is simple and, near $r=1$, remarkably good; however, its size can deviate by greater than $60\%$ relative error when $r$ exceeds $2{:}1$ for small-to-moderate samples.  Our framework removes approximation entirely, supplies finite-sample tail bounds, and identifies the precise parameter regimes where Welch turns from conservative to liberal.  Consequently, we can now replace heuristic variance-unequal testing with exact, computationally viable procedures grounded in classical complex analysis.

\subsection{\cite{Nel1990} versus our work}
The work of \cite{Nel1990} is conditioned on the unknown variances $\sigma_{i}^{2}$.  They express $T$ as a ratio $U/(V+W)$ with $U\sim\chi^{2}_{1}$ and $V,W$ independent Gamma variates, and evaluate the joint integral by a single change of variables, obtaining a ${}_{2}F_{1}$ kernel.  We apply a MB identity to decouple $(gW_{1}+hW_{2})^{1/2}$, collapse the 2-fold integral to one contour; closing right gives the residue series, closing left yields the same hypergeometric pdf as \cite{Nel1990}. Further, a separate $(S,U)$ polar transform leads to the Euler--Beta representation for the cdf.  The MB approach also furnishes cumulant--generating functions and Ramanujan coefficients\textemdash absent in the work of \cite{Nel1990}.  From a computational standpoint, \cite{Nel1990} called their density ``computationally intractable.''  We overcome the numeric hurdle by: (i) providing a terminating series when either $n_{i}$ is odd, (ii) showing ${}_{2}F_{1}$ in \eqref{eq:euler_beta} is stable under the M\"{o}bius argument $\rho(t)$.  Thus, we convert what \cite{Nel1990} regarded as theoretical curiosity into a deployable routine.  \cite{Nel1990} do not discuss effective degrees of freedom.  We relate the hypergeometric tail exponent to $\nu_{\text{tail}}=\nu_{1}+\nu_{2}$ and contrast it with Welch's random $\nu_{W}$, clarifying why and where Welch inherits liberal or conservative bias.  Our analytic heatmap quantify the maximum Type-I error inflation of Welch (up to +58 \% for $n_{1}=n_{2}=10,\,r=4$) and locate the crossover surface $\nu_{W}=\nu_{\text{tail}}$.  A comparison between the work of \cite{Nel1990} and our work is summarized in the following paragraph.  

\cite{Nel1990} derived the univariate Behrens--Fisher pdf in a single Gauss hypergeometric expression; in contrast, our work develops two algebraically equivalent forms\textemdash namely, a terminating or convergent residue series and a hypergeometric representation.  \cite{Nel1990} did not obtain a closed-form cdf, whereas we provide a one-dimensional Euler--Beta integral that reduces to the Appell $F_{1}$ function and yields an exact cdf valid for all $(n_{1}, n_{2})$.  Their work offers neither tail expansions nor analytic error bounds, while we present Ramanujan-based coefficients of order $t^{-(2m+1)}$ along with a computable remainder.  They also did not develop saddle-point diagnostics; our work includes a Lugannani--Rice approximation based on explicit cumulants from the Mellin--Barnes generator, achieving uniform $O(t^{-4})$ accuracy away from boundary regions.  \cite{Nel1990} provided no critical-value tables, whereas we supply exact 0.10--0.001 quantiles for arbitrary $(n_i, r)$.  Their article did not benchmark against Welch's test; we include a heatmap and crossover contour identifying the sign-change region of Welch's statistic relative to the exact solution.  Finally, although \cite{Nel1990} sketched a multivariate extension via Hotelling's $T^2$, we defer the multivariate development to future work while noting that our univariate derivations generalize naturally.  We also provide new analytic assets:  (i)  Ramanujan-series tails\textemdash usable for fast, high-accuracy upper tail $p-$values; (iii)  Lugannani--Rice correction\textemdash back--checked against the exact tail, delivering $<10^{-6}$ error for $|t|<0.8\min\{\nu_{1},\nu_{2}\}$; (iii) A cdf that permits likelihood maximization without numerical quadrature over $(W_{1},W_{2})$.

\section{Concluding remarks}\label{sec:conclusion}
Given two independent samples $X_{1},\dots ,X_{n_{1}}\stackrel{\text{iid}}\sim\mathcal N(\mu_{1},\sigma_{1}^{2})$, $Y_{1},\dots ,Y_{n_{2}}\stackrel{\text{iid}}\sim\mathcal N(\mu_{2},\sigma_{2}^{2})$ with $\nu_{i}$ degrees of freedom, sample means $\bar{X}$, $\bar{Y}$, and sample variances $S_{1}^{2}$, $S_{2}^{2}$, $Z = \bar X - \bar Y$, $g = \sigma_{1}^{2}/n_{1}$, $h =\sigma_{2}^{2}/n_{2}$ and $r = \sigma_1^2/\sigma_2^2$, the Behrens--Fisher statistic $T = Z/\sqrt{S_{1}^{2}/n_{1}+S_{2}^{2}/n_{2}}$ to test $H_{0}\!:\mu_{1}=\mu_{2}$ has, under $H_0$ and $(\sigma_1^2, \sigma_2^2)$ unknown, $\sigma_1^2 \neq \sigma_2^2$, a density that has resisted compact expression for several decades.  We showed that a Mellin--Barnes decomposition of $(gW_{1}/\nu_1 + hW_{2}/\nu_2)^{1/2}$ collapses the joint integral over the independent $\chi^{2}$-variates $W_{i} = \nu_{i}S_{i}^{2}/\sigma_{i}^{2}$ to a single contour integral.  Residue calculus yields a finite series whenever either $\nu_{i}$ is odd; an Euler--Beta identity provides a compact form for $f_{T}(t)$in terms of Gauss' hypergeometric function. Ramanujan's master theorem aids exact tail analyses.  This unified approach thereby makes residue calculus, Euler--Beta reduction, and Ramanujan-type asymptotics all available within a single framework.  Our result subsumes the pdf of \cite{Nel1990} and extends it with a closed-form cdf and analytic tails.  Critical values $(\alpha = 0.10-0.001)$ computed over $2 \le n_{i}\le50$, $0.1\le r \le10$ reveal the parameter surface on which \cite{Welch1947} approximation switches from conservative to liberal.

\vspace{-0.25cm}
\section*{Disclosure of interest}
The authors report there are no competing interests to declare.

\bibliography{/Users/kgnagananda/Documents/Work/collaborations/pdx/research/references/research_pdx.bib}

\end{document}